\documentclass[a4paper, reqno]{amsart}
\usepackage{a4,wasysym}
\usepackage{amsmath,graphicx, amssymb,color,ulem,bbm}

\newtheorem{theorem}{Theorem}
\newtheorem{corollary}[theorem]{Corollary}
\newtheorem{lemma}[theorem]{Lemma}

\newcommand{\ud}{\mathrm{d}}

\newcommand{\eofproof}{\hfill$\Box$ }
\begin{document}

\title[Structured populations with distributed recruitment]{Structured populations with distributed recruitment: from PDE to delay formulation}

\author[\`{A}. Calsina]{\`{A}ngel Calsina}
\address{\`{A}ngel Calsina, Department of Mathematics, Universitat Aut\`{o}noma de Barcelona, Bellaterra, 08193, Spain}
\email{acalsina@mat.uab.es}

\author[O. Diekmann]{Odo Diekmann}
\address{Odo Diekmann, Department of Mathematics, University of Utrecht, Budapestlaan 6, PO Box 80010, 3508
TA, Utrecht, The Netherlands}
\email{O.Diekmann@uu.nl}

\author[J. Z. Farkas]{J\'{o}zsef Z. Farkas}
\address{J\'{o}zsef Z. Farkas, Division of Computing Science and Mathematics, University of Stirling, Stirling, FK9 4LA, United Kingdom }
\email{jozsef.farkas@stir.ac.uk}

\date{\today}

\begin{abstract}
In this work first we consider a physiologically structured population model with a distributed recruitment process. That is, our model allows newly recruited individuals to enter the population at all possible individual states, in principle. The model can be naturally formulated as a first order partial integro-differential equation, and it has been studied extensively. In particular, it is well-posed on the biologically relevant state space of Lebesgue integrable functions. We also formulate a delayed integral equation (renewal equation) for the distributed birth rate of the po\-pu\-lation. We aim to illustrate the connection between the partial integro-differential and the delayed integral equation formulation of the model utilising a recent spectral theoretic result. In particular, we consider the equivalence of the steady state problems in the two different formulations, which then leads us to characterise irreducibility of the semigroup governing the linear partial integro-differential equation. Furthermore, using the method of characteristics, we investigate the connection between the time dependent problems. In particular, we prove that any (non-negative) solution of the delayed integral equation determines a (non-negative) solution of the partial differential equation and vice versa. The results obtained for the particular distributed states at birth model then lead us to present some very general results, which establish the equivalence between a general class of partial differential  and delay equation, modelling physiologically structured populations.
\end{abstract}
\maketitle

\section{Prologue}

Structured population models are of great interest both from the mathematical and the application point of view.
Traditionally they have been formulated as partial differential equations (PDEs for short) often with non-local (and nonlinear) boundary conditions. The early monograph \cite{W} provided a comprehensive mathematical treatment of nonlinear age-structured population dynamics; while \cite{CUS} provides valuable insight into the modelling and analysis of size-structured populations, which is accessible for theoretical biologists alike. We also refer the interested reader to the ``green book'' \cite{MD}, which provides a systematic introduction to structured population models. 

Today we can safely say that the theory of semilinear equations, i.e. age-structured models is well understood. The monograph \cite{W} by G. F. Webb  provides a general ma\-the\-matical treatment of a wide class of models; see also the very useful  monograph \cite{I}. In \cite{W}, among other things, the Principle of Linearised Stability (PLS for short) was established for a general class of age-structured (semilinear) models. We also refer the interested reader to the paper \cite{Grabosch}, where generators of translation semigroups were cha\-rac\-te\-rised. The general results obtained in \cite{Grabosch} can be applied directly to a number of structured population models. 

On the other hand, the analysis of quasilinear models, which naturally arise when modelling size-structured populations, still poses difficult  mathematical challenges. Remarkably, the PLS has not been established yet for fully general quasilinear PDE models. Partly because of this shortfall, a number of researchers have worked on establishing a general theory using a delay/integral equation formulation of physiologically structured population models, see e.g. \cite{D2,D1,DGM2,DGM1}, and the references therein. As a result, a comprehensive modelling approach to build physiologically structured population models from basic principles was developed. For a class of these models, applying the sun-star calculus for abstract integral equations, the Principle of Linearised Stability and the Hopf-bifurcation theorem were established in \cite{D2,D1}. For a comprehensive overview of the underlying mathematical theory see e.g. \cite{DGLW}. We also mention, that an alternative approach uses integ\-ra\-ted semigroups, see \cite{Magal}, and the references therein.

It is fair to say that both (delay and PDE) modelling approaches can be fruitfully applied to analyse concrete problems. The reader may want to compare for example the papers \cite{DGMNR} and \cite{FH} for two different formulations (and analysis) of a classic Daphnia-algae (structured consumer-resource) model. Nevertheless, to our knowledge, the literature lacks rigorous results establishing the connection between the two modelling approaches. One may anticipate that, at least in a simple setting, the connection is quite easily understandable and can be motivated naturally, either from the mathematical, or from the biological point of view. The ever increasing volume of the literature on structured po\-pu\-lation dynamics also underlines the need of a better understanding of the connection between the two modelling approaches. 

In case of a classic example, the basic linear age-structured model formulated by McKendrick in the early 20th century, the integral equation formulation arises naturally for example by solving the partial differential equation using the method of characteristics. The method of characteristics, to recast the partial differential equation as a pair of in\-teg\-ral equations, can be applied to more general (including nonlinear) models, too; but the mathematical machinery often lacks a biologically inspired motivation.

In this work we are going to consider models of physiologically structured populations with a distributed recruitment process, 
that is, we allow a continuum of states at birth. A vast number of models with distributed recruitment processes have been introduced mainly as models of cell populations, we just mention here some of the early works, e.g. \cite{Heijmans, MD, WebbGrabosch}, for the interested reader. More recently, we have been investigating size-structured models with a distributed recruitment process. In the papers \cite{AF,BCC,CF2,CF,FGH,FH} the emphasis was on qualitative questions, such as existence and stability of steady states; while in \cite{CS2} well-posedness of a general class of models was established. In this work our  main goal is to establish a rigorous connection between solutions of the partial differential equation and of the delayed integral equation formulation of the distributed states at birth model. Perhaps somewhat interestingly, we are going to show how a fairly recent spectral theoretic result may help us to understand the connection (at least formally) between the two formulations. With this, we also hope to ignite researchers to investigate the connection between the two different formulations, for different classes of structured population models. Any progress in this direction may not just contribute to the general mathematical theory of structured population dynamics, but it may also shed light on some interesting properties of the models. We just mention here as an example that net reproduction numbers/functions/functionals naturally arise when considering steady state problems (and addressing stability questions) in a PDE formulation, see e.g. \cite{FH2,FH}. At the same time the net reproduction number, often denoted by $R_0$, also arises in the cumulative formulation as the spectral bound of the next generation operator, see e.g. \cite{DGM,Inaba} for more details. 

\section{Motivation and speculations}

First we are going to study a prototype linear model of a physiologically structured population with distributed recruitment, i.e. in our model individuals may enter the population at all possible individual states. We refer the reader to \cite{CS2}, where a global existence result for a very general nonlinear model with distributed recruitment process was established. Qualitative questions, such as existence and stability of steady states, and asynchronous exponential growth, were addressed in \cite{AF,FGH}; while in \cite{AFLM} a finite difference scheme was developed for a nonlinear distributed states at birth model to aid numerical simulations.

Our motivation stems from the very recent papers \cite{AF,FGH}, in which we considered the positive steady state problem
for a nonlinear physiologically structured population model with distributed recruitment process. In the earlier paper \cite{FGH}, following ideas from \cite{BCC}, we treated the steady state problem of the model using semigroup methods, via analysing spectral properties of the unbounded semigroup generator. In contrast, in \cite{AF}, we reformulated and studied the steady state problem in the form of a parametrised fa\-mi\-ly of integral equations. This then led us to arrive at the definition of a particular net reproduction function (motivated by the mathematical formulation). We then investigated how this net reproduction function, which is naturally related to the existence of positive steady states of the nonlinear model, is connected to a biologically meaningful and relevant net reproduction function. In the light of the results obtained in \cite{AF,FGH}, it seems to be a very natural and interesting question to consider whether the two formulations of the positive steady state problem are equivalent. Furthermore, the question naturally arises, whether one can establish the equivalence of the corresponding time-dependent problems (PDE and delayed integral equation). We are going to answer these questions in this work. With the aim of presenting results accessible to the widest possible audiences, we start with a specific model and recall some of the important ideas presented in \cite{AF}. 

The following partial differential equation may be considered as the basic physiologically structured population model with distributed recruitment.
\begin{equation}\label{distributed}
\begin{aligned}
p_t(s,t)+\left(\gamma(s)p(s,t)\right)_s &= -\mu(s)p(s,t)+\int_0^1\beta(s,y)p(y,t)\,\ud y,\quad s\in (0,1), \\
\gamma(0)p(0,t)&=0.
\end{aligned}
\end{equation}
Our model \eqref{distributed} is equipped with an initial condition $p(\cdot,0)$, which determines the initial population density. In this work, for the ease of presentation, we use the interval $[0,1]$ for the structuring variable $s$.  
The model above describes the time evolution of the size-distribution, (or the distribution of any another physiological structuring variable), of a population with distributed recruitment process. That is, individuals may be recruited into the population at different sizes/states. More precisely, the recruitment rate is determined  by the fertility function $\beta$. As usual, $\gamma$ and $\mu$ denote individual growth and mortality rates, respectively. All of the vital rates depend on the structuring variable $s$, in general.

We impose some natural regularity assumptions on the model ingredients: $\beta,\mu,\gamma$. Specifically, we require $\gamma$ to be strictly positive, continuously differentiable, and we assume that $\mu$ and $\beta$ are non-negative and continuous. 
These conditions suffice to analyse the problem both in the framework of partial differential and integral equation. We note that the assumption of strict positivity on $\gamma$ implies that in principle there may be a number of individuals growing beyond the maximal size. However, choosing the mortality function such that it becomes sufficiently large close to the maximal size implies that the proportion of individuals growing beyond the maximal size is negligible. At the same time we note that our aim here is not to handle the most general situation (i.e. to study the model with the most relaxed smoothness assumptions on the vital rates).

For problem \eqref{distributed} the natural choice of state space is the Lebesgue space $L^1(0,1)$. 
This is simply because convergence in the $L^1$ norm has a clear biological interpretation, that is, it implies convergence of the density describing the size distribution of the population. 
On this state space, under the assumptions on $\beta,\mu$ and $\gamma$ we specified earlier, the following densely defined operator generates a strongly continuous semigroup of bounded linear operators.
\begin{align}\label{Aopdef}
\mathcal{A}\, p=-\frac{\partial}{\partial s}\left(\gamma p\right)-\mu p+\int_0^1\beta(\cdot,y)p(y)\,\ud y,
\end{align}
where the domain of the generator $\mathcal{A}$ is defined as 
\begin{equation*}
D(\mathcal{A})=\left\{p\in W^{1,1}(0,1)\,|\, \gamma(0)p(0)=0\right\}.
\end{equation*}
It is shown (see e.g. \cite{FGH}) that $\mathcal{A}$ generates an eventually compact and positive semigroup.
Moreover, it was shown in \cite{FGH}, that the semigroup is irreducible, if there exists an $\varepsilon_*>0$ such that $\forall\,\varepsilon\in (0,\varepsilon_*)$ we have 
\begin{align}
\int_0^{\varepsilon}\int_{1-\varepsilon}^1\beta(s,y)\,\ud y\,\ud s>0.\label{irredcond}
\end{align}
Note that this irreducibility condition has a clear biological meaning: it requires that large individuals produce offspring of arbitrarily small size. For example if $\beta$ is continuous (as indeed we assumed), condition \eqref{irredcond} holds if $\beta(0,1)\ne 0$. Also note that condition \eqref{irredcond} is sufficient, but not necessary. For example the weaker assumption, that there exists an $s_*\in [0,1]$ and an $\hat{s}\in [s_*,1]$, such that $\beta(s_*,1)>0$ and $\beta(0,\hat{s})>0$ hold simultaneously, would also suffice to guarantee irreducibility of the semigroup. In fact, we are going to return to the question of irreducibility of the semigroup, and we will formulate a necessary and sufficient condition in the next section.

Due to the eventual compactness of the governing semigroup, the spectrum of $\mathcal{A}$ contains only eigenvalues of finite multiplicity, and the Spectral Mapping Theorem holds true. Moreover, the eigenvalues of $\mathcal{A}$ can be characterised implicitly via the following functional equation
\begin{equation}\label{eigeninteq}
f(\cdot)=\displaystyle\int_0^1\frac{\beta(\cdot,y)}{\gamma(y)}\int_0^yf(r)\exp\left\{-\int_r^y\frac{\lambda+\mu(x)}{\gamma(x)}\,\ud x\right\}\,\ud r\,\ud y,
\end{equation}
where we introduced
\begin{equation*}
f(\cdot)=\int_0^1\beta(\cdot,y)p(y)\,\ud y.
\end{equation*}
More precisely, $\lambda\in\mathbb{C}$ is an eigenvalue of $\mathcal{A}$, if and only if there exists a non-trivial $f$, such that $\left(\lambda,f\right)$ is a solution of equation \eqref{eigeninteq}. Then, the non-trivial eigenvector $p$ corresponding to $\lambda$, is given by
\begin{equation}\label{eigeninteq2}
p(\cdot)=\frac{1}{\gamma(\cdot)}\int_0^\cdot\exp\left\{-\int_r^\cdot\frac{\lambda+\mu(x)}{\gamma(x)}\,\ud x\right\}f(r)\,\ud r.
\end{equation}
The sign of the spectral bound of $\mathcal{A}$, (note that the spectral bound is a real and strictly dominant eigenvalue), determines the asymptotic behaviour of solutions of the linear model \eqref{distributed}. 
Note that equation \eqref{eigeninteq}, which is the characteristic equation corresponding to the linear problem \eqref{distributed}, is a functional equation, in contrast to the case of a single state at birth model, when the characteristic equation is a scalar one. 

We also note that even in the case of an unbounded size-space, the asymptotic behaviour of solutions of the linear model (or a linearised model),  is determined by eigenvalues at least if the mortality function is bounded below by a positive constant $\mu_0$. In this case, even though the essential spectrum of $\mathcal{A}$ may not be empty, it can still be contained in the half-plane $\lambda\in\mathbb{C}$, Re$(\lambda)\le -\mu_0$. We refer the interested reader for example to \cite{FH3} for a result of this kind, established for a size-structured population model incorporating cannibalistic behaviour. 

Equation \eqref{eigeninteq} naturally leads us to define a bounded integral operator $\mathcal{L}$ as follows. 
\begin{equation}
\mathcal{L}\, b=\int_0^1\frac{\beta(\cdot,y)}{\gamma(y)}\int_0^yb(r)\exp\left\{-\int_r^y\frac{\mu(x)}{\gamma(x)}\,\ud x\right\}\,\ud r\,\ud y,\label{intop1}
\end{equation}
with domain $L^1(0,1)$. Equivalently, by introducing the function $b$ as
\begin{equation*}
b(\cdot)=\int_0^1\beta(\cdot,y)p(y)\,\ud y,
\end{equation*} 
and solving the equation $\mathcal{A}\,p=0$, one arrives at the integral equation $\mathcal{L}\, b=b$, with $\mathcal{L}$ defined as above in \eqref{intop1}. It is then shown that the operator $\mathcal{A}$ has spectral bound $0$, (which is a simple eigenvalue if the semigroup is irreducible, e.g. when \eqref{irredcond} holds), if and only if the integral operator $\mathcal{L}$ has spectral radius $1$. This result is in the spirit of \cite[Ch.5]{MD} by Heijmans, but can be obtained directly by using a fairly recent spectral theorem from \cite{HT} due to Horst Thieme. We recall now this result for the reader's convenience. In the theorem below $s$ stands for the spectral bound of a linear operator, while $r$ denotes the spectral radius of a (bounded) linear operator. An operator $\mathcal{O}$ is called resolvent-positive if its resolvent set contains an interval of the form $(\omega,\infty)$, and the operator $(\lambda\mathcal{I}-\mathcal{O})^{-1}$ is positive for $\lambda$ large enough.
\begin{theorem}\label{spectralth} \cite[Theorem 3.5]{HT}
Let $\mathcal{B}$ be a resolvent-positive operator on $\mathcal{X}$, $s(\mathcal{B})<0$, and $\mathcal{A} =\mathcal{B}+\mathcal{C}$ a positive perturbation of $\mathcal{B}$, (i.e. $\mathcal{C}$ is a positive linear operator with domain including the domain of $\mathcal{B}$). If $\mathcal{A}$ is resolvent-positive then $s(\mathcal{A})$ has the same sign as $r\left(-\mathcal{C}\,\mathcal{B}^{-1}\right)-1$.
\end{theorem}
To apply Theorem \ref{spectralth} we split the operator $\mathcal{A}$ defined in \eqref{Aopdef} into two parts. Hence, we rewrite the partial differential equation \eqref{distributed} as the following abstract Cauchy problem.
\begin{equation}\label{Cauchy-abstr}
\frac{\ud p}{\ud t}=\mathcal{A}\,p=\left(\mathcal{B}+\mathcal{C}\right)p, \quad
p(0)=p^0,
\end{equation}
where we define
\begin{align}
\mathcal{B}\,p&=-\frac{\partial}{\partial s}\left(\gamma\,p\right)-\mu p,\quad D(\mathcal{B})=D(\mathcal{A})=\left\{p\in W^{1,1}(0,1)\,|\, \gamma(0)p(0)=0\right\},\label{Bdef} \\
\mathcal{C}\,p&=\int_0^1\beta(\cdot,y)p(y)\,\ud y,\quad D(\mathcal{C})=L^1(0,1).\label{Cdef}
\end{align}
Note that $\mathcal{C}\,:\,D(\mathcal{B})\to\mathcal{X}$ is positive in the sense as defined in \cite{HT}. That is, for any 
$p\in D(\mathcal{B})\cap\mathcal{X}_+$ we have that $\mathcal{C}\, p\in\mathcal{X}_+$. Then, $-\mathcal{C}\mathcal{B}^{-1}$ is the integral operator $\mathcal{L}$ defined in \eqref{intop1}. This was shown in \cite{AF}. If the mortality $\mu$ is not identically zero, then it can be shown that the spectral bound of $\mathcal{B}$ is negative, hence the operators $\mathcal{A,B}$ and $\mathcal{C}$ satisfy the conditions of Theorem \ref{spectralth}.
That is, the connection between the two eigenvalue problems
\begin{equation}\label{eigvss}
\mathcal{A}\,p=0\,p,\quad \text{and}\quad \mathcal{L}\,b=1\,b,
\end{equation}
is given by Theorem \ref{spectralth}. In principle, any positive eigenvector of $\mathcal{A}$ corresponding to the eigenvalue $0$ determines a positive steady state of the PDE \eqref{distributed}. Similarly, any positive eigenvector $b$ of $\mathcal{L}$ corresponding to the eigenvalue $1$ will be a steady state of the delayed integral equation, see also \eqref{inteq1} below.
At the same time, any positive eigenvector $b$ of $\mathcal{L}$ corresponding to the eigenvalue $1$, determines a positive steady state $p$ 
of the PDE \eqref{distributed} as 
\begin{equation}\label{steadyconn}
p(\cdot)=\frac{1}{\gamma(\cdot)}\int_0^\cdot b(r)\exp\left\{-\int_r^\cdot\frac{\mu(x)}{\gamma(x)}\,\ud x\right\}\,\ud r.
\end{equation}
In the next section we will further investigate whether we can conclude that the equi\-va\-lence of the two eigenvalue problems implies that the (positive) steady state problems in the two different formulations are equivalent. We note however, that studying the eigenvalue problems \eqref{eigvss} does not allow us to give the best possible answer for this question, as we will see later in Section 5. However, this method is applicable to study existence of steady states of nonlinear models, see for example \cite{AF,FGH}; while we anticipate that the extension of the general results of Section 5 for nonlinear models will present considerable challenges. 

Our main goal in the rest of the section is to illustrate how Theorem \ref{spectralth} can lead our intuition to recover from the PDE model \eqref{distributed} the delayed integral equation formulation of the model. Of course the problems are time-dependent, and naturally, the delay equation is formulated on a different state space. This difference is very much the same though as for the time-independent problems \eqref{eigvss}. For related developments we also refer the interested reader to the earlier papers \cite{DGM,DGM1}. 

We start by deducing the delay formulation from basic modelling principles. Similarly to the case of a classic age or size-structured model with single state at birth, it is natural to formulate an integral equation for the following recruitment (or birth) function
\begin{equation}\label{birthrate}
b(s,t)=\int_0^1\beta(s,y)p(y,t)\,\ud y,\quad s\in [0,1],\,\,t>0.
\end{equation}
Since $\gamma$ denotes individual growth rate, i.e. the rate of change of $s$ with respect to time $t$, we can compute the
time $\tau(x,y)$ an individual spends to grow from size $x$ to size $y$. We have
\begin{align*}
\frac{\ud s}{\ud t}=\gamma(s(t)),\quad s(0)=x,
\end{align*}
which yields
\begin{align*}
\int_x^y\frac{1}{\gamma(r)}\,\ud r=\int_0^\tau 1\,\ud t=\tau(x,y).
\end{align*}
On the other hand, the probability to survive from size $x$ to size $y$ is given by
\begin{align*}
\exp\left\{-\int_x^y\frac{\mu(r)}{\gamma(r)}\,\ud r\right\}.
\end{align*}
Hence the flux of individuals of size $y$ at time $t$ is
\begin{align*}
\gamma(y)p(y,t)=\displaystyle\int_0^yb\left(x,t-\int_x^y\frac{1}{\gamma(r)}\,\ud r\right)\exp\left\{-\int_x^y\frac{\mu(r)}{\gamma(r)}\,\ud r\right\}\,\ud x.
\end{align*}
The formula above is meaningful at least for times $t>\displaystyle\int_0^1\frac{1}{\gamma(s)}\,\ud s=:\Gamma$. This is because in principle the density $p$, which defines $b$ in \eqref{birthrate} is not defined for negative times. Hence from \eqref{birthrate} we obtain that the birth rate function $b$ satisfies the integral equation (now with delay)
\begin{align}
b(\cdot,t)=\int_0^1\frac{\beta(\cdot,y)}{\gamma(y)}\int_0^yb\left(x,t-\int_x^y\frac{1}{\gamma(r)}\,\ud r\right)
\exp\left\{-\int_x^y\frac{\mu(r)}{\gamma(r)}\,\ud r\right\}\,\ud x\,\ud y, \label{inteq1}
\end{align}
for $t>\Gamma$.

We now intend to illustrate how the integral equation \eqref{inteq1} can be obtained (formally) in the spirit of Theorem \ref{spectralth} by introducing some appropriate operators. Note that what follows will be done by hand-waving (on purpose, to motivate the connection to Theorem \ref{spectralth}), and the rigorous connection between the time-dependent problems will be given in  Section 4. For a given initial condition $p^0$, let $p$ be the solution of the PDE \eqref{distributed}. It is natural to assume that $p\in C\left([0,\infty);L^1(0,1)\right)$. First we obtain the recruitment (as in \eqref{birthrate}) function from the population density distribution now formally, by defining a map $C$ as follows.
\begin{align}
C \,:& \,\underbrace{p(\diamond,\ast)}_{\in C\left([0,\infty);L^1(0,1)\right)} \longrightarrow\underbrace{b(\cdot,\ast)}_{\in C\left([0,\infty);C(0,1)\right)}\left(=\int_0^1\beta(\cdot,\diamond)p(\diamond,\ast)\,\ud \diamond\right).
\end{align}
Note that above $\diamond$ stands for the size variable (first component of $p$) and $*$ stands for the time variable (the second component of $p$). At the same time (as noted above) $p$ can be understood as a continuous $L^1$-valued function of time. Formally, $C$ can be considered to be the time dependent analogue of the linear operator $\mathcal{C}$ defined in \eqref{Cdef}.  Next we define a shift operator, which is required, in contrast to the steady state problem, due to the time-dependence. In particular, we introduce
\begin{align}
S\,:& \, \underbrace{b(\cdot,\ast)}_{\in C\left([0,\infty);C(0,1)\right)} \longrightarrow\quad\underbrace{b\left(\cdot,\ast-\int_\cdot^\diamond\frac{1}{\gamma(r)}\,\ud r\right)}_{\in C\left([\Gamma,\infty);C([0,1]^2)\right)}.
\end{align}
Similarly, we define a time-dependent version of the operator $-\mathcal{B}^{-1}$ as follows.
\begin{align}
-B^{-1}\,:&\, \underbrace{b(\cdot,\ast,\diamond)}_{\in C\left([\Gamma,\infty);C([0,1]^2)\right)} \longrightarrow\quad\underbrace{-\frac{1}{\gamma(\diamond)}\int_0^\diamond\exp\left\{-\int_\cdot^\diamond\frac{\mu(r)}{\gamma(r)}\,\ud r\right\}b(\cdot,\ast,\diamond)\,\ud \cdot}_{\in C\left([\Gamma,\infty);L^1(0,1)\right)}.
\end{align}
Then, the delayed integral equation \eqref{inteq1} can formally (and in an abstract sense) be written as
\begin{align}
b=\left(-CB^{-1}S\right)\,b.
\end{align}

Note that it is straightforward to verify that of course the integral equation \eqref{inteq1} is indeed the one 
which arises when building the linear distributed states model following the steps described by Diekmann et al., see e.g. \cite{D2,DGLW,DGM1}. 

We would also like to point out that one would readily obtain the following (similar) delayed integral equation for the birth rate 
\begin{align}
b(\cdot,t)=&\displaystyle\int_0^1\frac{\beta(\cdot,y,P(t))}{\gamma(y)} \nonumber \\
& \quad\times\displaystyle\int_0^yb\left(s,t-\int_s^t\frac{1}{\gamma(r)}\,\ud r\right)
\exp\left\{-\displaystyle\int_s^y\frac{\mu\left(r,P\left(t-\int_s^r\frac{1}{\gamma(x)}\,\ud x\right)\right)}{\gamma(r)}\,\ud r\right\}\,\ud s\,\ud y, \label{inteq2}
\end{align}
for $t>\Gamma$, in the semilinear case, i.e. when $\mu$ and $\beta$ depend on the total population size, $P(t)=\int_0^1 p(s,t)\,\ud s$, too. 
The quasilinear case, i.e. when $\gamma$ depends on the total population size as well, leads to a formally more complicated (and cumbersome) delayed integral equation. Note that, even in the case when $\gamma$ depends on the total population size $P$, the partial differential equation can still be solved along characteristic lines (at least locally), see e.g. \cite{Calsina,G-H,Kato}.

\section{Exploring the equivalence of the steady state problems}

As we have seen in the previous section, the eigenvalue problems \eqref{eigvss} are very closely related by Theorem \ref{spectralth}. Here we further investigate the connection between these two eigenvalue problems, and we also discuss the relationship between the two steady state problems. To this end, note that if the semigroup generated by $\mathcal{A}$ is irreducible, then the spectral bound of $\mathcal{A}$ is the only eigenvalue with a positive eigenvector, see e.g. Proposition 3.5 and Theorem 3.8 in \cite[C-III]{AGG}, as well as Theorem 1.2 and Remark 1.3 in \cite{AB}. Moreover, the (unique, normalised) eigenvector corresponding to the spectral bound is strictly positive. Similarly, if the integral operator $\mathcal{L}$ is irreducible, then the spectral radius of $\mathcal{L}$ is the only eigenvalue with a corresponding strictly positive eigenvector, see e.g. Theorem 5.2 and its corollary  in \cite[Ch.V]{Sch}. Hence if the irreducibility of the semigroup would imply irreducibility of the integral operator defined in \eqref{intop1} and vice versa, we could say that the two spectral problems \eqref{eigvss} are equivalent. Therefore, it is very interesting to explore what is the relationship between the irreducibility of the semigroup generated by $\mathcal{A}$, and the irreducibility of the integral operator $\mathcal{L}$. Note that irreducibility of the semigroup generated by $\mathcal{A}$ plays an important role in the qualitative analysis of the nonlinear version of model \eqref{distributed}, too, see e.g.  \cite{FGH}. At the same time we note that the existence of a unique strictly positive eigenvector corresponding to the eigenvalue (spectral bound) $0$ of $\mathcal{A}$ is a sufficient, but it is clearly not a necessary condition, for the existence of a positive steady state of problem \eqref{distributed}. 

We first formulate a necessary and sufficient condition for the irreducibility of the semigroup generated by $\mathcal{A}$, for a continuous and non-negative $\beta$. Let us recall from \cite[C-III]{AGG} that a positive semigroup $\mathcal{T}$ on the Banach lattice $L^1(0,1)$ is said to be irreducible, if for any given $f\in L_+^1(0,1)$ and $\phi\in L_+^\infty(0,1)$, we have $\langle \mathcal{T}(t_*)\,f,\phi\rangle >0$ for some $t_*\ge 0$, where $\langle \cdot,\cdot\rangle$ stands for the natural pairing between $L^1$ and its dual $L^\infty$. 

We also recall another equivalent characterisation of irreducibility (see again \cite[C-III]{AGG}). In particular, the semigroup $\mathcal{T}$ is irreducible if it does not admit invariant closed ideals other than the trivial ones. We also note that in $L^1(0,1)$ the non-trivial ideals are characterised by equivalence classes of functions vanishing on a measurable subset of $(0,1)$ with positive measure, see e.g. \cite{AGG}. Our main result is formulated in the following theorem. 
\begin{theorem}
The semigroup $\mathcal{T}$ generated by the operator $\mathcal{A}$ defined in \eqref{Aopdef} is irreducible, if 
and only if 
\begin{equation}
\forall \alpha\in (0,1),\,\, \exists\, (s_*,y_*)\in [0,\alpha]\times [\alpha,1],\quad  \text{such that}\quad \beta(s_*,y_*)>0.\label{irredchar}
\end{equation}
\end{theorem}
{\bf Proof.}
Note that is not too difficult to see that the first definition of irreducibility, which we recalled earlier, is equivalent to the following condition: for every $0\not\equiv p_0\in L_+^1(0,1)$ one has that
\begin{equation}
\displaystyle\bigcup_{t\ge 0} ess\,supp\left(\mathcal{T}(t)\, p_0\right)= [0,1].\label{irredeq1}
\end{equation}

First to see the necessity of condition \eqref{irredchar}, note that if there was an $\alpha\in (0,1)$, such that for every 
$(s,y)\in [0,\alpha]\times [\alpha,1]$, $\beta(s,y)=0$ would hold; then the ideal consisting of functions vanishing almost everywhere on the interval $(0,\alpha)$ would be invariant under the semigroup $\mathcal{T}$. 

To prove that condition \eqref{irredchar} is also sufficient, we note that for any initial condition $0\not\equiv p_0\in L_+^1 (0,1)$, if for some $\hat{t}\ge 0$, $s_*\in \text{ess supp}\left(\mathcal{T}(\hat{t})\,p_0\right)$, then $\forall\,s>s_*,\, \exists\, t\ge\hat{t}$, such that $s\in\,\text{ess supp}\left(\mathcal{T}(t)\, p_0\right)$. Note that, this in particular implies that for any $0\not\equiv p_0\in L_+^1(0,1)$ 
we have 
\begin{equation*}
\displaystyle\bigcup_{t\ge 0} ess\,supp\left(\mathcal{T}(t)\, p_0\right)= [s_{min},1],
\end{equation*}
for some $s_{min}\in [0,1)$. We are going to prove that in fact $s_{min}>0$ is not possible. To this end, we define a function $R$ as follows.
\begin{equation}\label{magicfunction}
R(x)=\inf\left\{s\in [0,1)\,|\, \beta(s,\tau)>0, \quad \text{for some}\quad \tau\ge x\right\},\quad x\in [0,1).
\end{equation}
First note that assumption \eqref{irredchar} and the continuity of $\beta$ imply that $R(0)=0$, and that $R(x)<x$ holds, for every $x\in (0,1)$. 

For any given initial condition $p_0\in L_+^1(0,1)$, starting with any $x_0\in (0,1)$ from the essential support of $p_0$, we define a sequence as follows.
\begin{equation*}
x_{n+1}:=R(x_n),\quad n\ge 0.
\end{equation*} 
Note that, $\forall n\ge 0$, $x_n\in\, \text{ess supp}\left(\mathcal{T}(t)\,p_0\right)$, for some $t\in [0,n\Gamma]$, hence 
$x_n\ge s_{min}, \, \forall\,n\in\mathbb{N}$. Since the sequence $x_n$ is decreasing (and bounded), it is convergent. The only fixed-point of $R$ is $0$. Hence, to show that  $x_n$ tends to $0$ we will prove that $R$ is right-continuous on $[0,1)$. 

Assume by the way of contradiction that $R$ is not right-continuous at some point $x_*\in [0,1)$. That is, $\exists\,\varepsilon>0$ such that $R\left(x_*+\frac{1}{n}\right)-R(x_*)>\varepsilon,\, \forall\, n\in\mathbb{N}$, and $n>n_*$ such that $x_*+\frac{1}{n}<1$. 
 
Then, it follows from the definition of $R$, that $\exists\, s_*\in \left(R(x_*),R\left(x_*+\frac{1}{n}\right)\right)$, and $\tau\ge x_*$ such that $\beta(s_*,\tau)>0$ 
holds, $\forall\,n>n_*$.
 
Since $\beta$ is continuous on the square $[0,1]\times [0,1]$, $\exists\,n'>n_*$ such that we have $\beta\left(s_*,\tau+\frac{1}{n'}\right)>0$.  This implies $R\left(x_*+\frac{1}{n'}\right)\le s_*$, a contradiction.

Therefore, since $x_n$ is decreasing, we have 
\begin{equation*}
\lim_{n\to \infty} x_{n+1}=\lim_{n\to\infty} R(x_n)=R\left(\lim_{n\to\infty} x_n\right),
\end{equation*}
that is, $\displaystyle\lim_{n\to\infty} x_n$ is the fixed point $0$ of $R$.  Therefore, it follows that $s_{min}=0$, 
hence \eqref{irredeq1} holds, and the proof is completed.
\eofproof

\begin{remark}
Note that the definition of the function $R$ in \eqref{magicfunction} can be clearly motivated from the biological point of view. Indeed, for any size $s\in (0,1)$, the value $R(s)$ is the smallest possible size of offspring produced by individuals of size $s$ in the rest of their lifetime.
\end{remark} 

Next we study the irreducibility of the integral operator $\mathcal{L}$ defined in \eqref{intop1}. First we note that it follows from \cite[Ch.V]{Sch} that $\mathcal{L}$ is irreducible, if and only if for every set $I\subset [0,1]$ of positive Lebesgue measure, we have
\begin{align}\label{Lirredcond}
\displaystyle\int\limits_{[0,1]\setminus I}\int\limits_{I} \left(\int_y^1\frac{\beta(s,r)}{\gamma(r)}\exp\left\{-\int_y^r\frac{\mu(\sigma)}{\gamma(\sigma)}\,\ud \sigma\right\}\,\ud r\right)\ud y\,\ud s>0.
\end{align}
If there was an $s_*\in (0,1)$, such that $\beta(s,y)=0,\, \forall\, (s,y)\in [0,s_*]\times [s_*,1]$, then with 
$I=[s_*,1]$, condition \eqref{Lirredcond} would not hold. Hence we conclude that criterion \eqref{irredchar} is necessary for the irreducibility of the integral operator $\mathcal{L}$. 

On the other hand, criterion \eqref{irredchar} is clearly not sufficient for the irreducibility of the integral operator $\mathcal{L}$. For example, any continuous and non-negative $\beta$, strictly positive on the open rectangle $(0,0.1)\times (0,1)$, but vanishing on $[0.1,1]\times [0,1]$, would satisfy condition \eqref{irredchar}, but clearly not condition  \eqref{Lirredcond}. 

On the other hand, since the kernel function in the double integral in \eqref{Lirredcond} is con\-ti\-nuous on the square $(0,1)\times (0,1)$, it is clear that the condition 
\begin{equation}\label{Lirredcond2}
\displaystyle\int_s^1\frac{\beta(s,r)}{\gamma(r)}\exp\left\{-\int_s^r\frac{\mu(\sigma)}{\gamma(\sigma)}\,\ud \sigma\right\}\,\ud r>0,\quad \forall\,s\in [0,1), 
\end{equation}
does imply condition \eqref{Lirredcond}. Moreover, since the function $\frac{1}{\gamma(\cdot)}\exp\left\{-\int_\diamond^{\cdot}\frac{\mu(\sigma)}{\gamma(\sigma)}\,\ud \sigma\right\}$ is strictly positive, condition \eqref{Lirredcond2} is equivalent to
\begin{equation}\label{Lirredcond3}
\int_s^1\beta(s,r)\,\ud r>0,\quad \forall\,s\in [0,1).
\end{equation}
This (sufficient) condition of irreducibility of $\mathcal{L}$ also has a clear biological interpretation. It requires that offspring of any size $s$ is "produced" by individuals of some larger size. 

We summarize formally our findings on irreducibility in the following corollary.
\begin{corollary}
Irreducibility of the integral operator $\mathcal{L}$ implies irreducibility of the semigroup $\mathcal{T}$ generated by the operator $\mathcal{A}$, but not vice versa.
\end{corollary}
In the light of the results above, let us give (for now) a partial answer for the question of equivalence of the PDE formulation \eqref{distributed}, and the delayed integral equation formulation \eqref{inteq1} (see also \eqref{de} below), at a steady state.  
 \begin{enumerate}
\item If both the semigroup generated by $\mathcal{A}$ and the integral operator $\mathcal{L}$ are irreducible, for example if condition \eqref{Lirredcond3} holds, then there are unique strictly positive (normalised) vectors $p$ and $b$ satisfying \eqref{eigvss}, and by Theorem \ref{spectralth} the steady state problems in the two different formulations are equivalent.  
\item If the semigroup generated by $\mathcal{A}$ is irreducible, but $\mathcal{L}$ is not irreducible, then there is still a unique strictly positive (normalised) vector $p$ satisfying the first equation in \eqref{eigvss}. This $p$ then determines a unique non-negative (but perhaps not strictly positive) $b$ satisfying the second equation in \eqref{eigvss} as 
\begin{equation}\label{steadyconn2}
b(\cdot)=\int_0^1\beta(\cdot,y)p(y)\,\ud y.
\end{equation}
In fact formula \eqref{steadyconn} shows that even though $\mathcal{L}$ may not be irreducible it does not have another non-negative eigenvector 
corresponding to its spectral radius $1$. Hence once again we may conclude that the steady state problems are equivalent. 
\item In the case when neither the semigroup generated by $\mathcal{A}$, nor the integral operator $\mathcal{L}$ are irreducible, the geometric multiplicities of the spectral bound $0$ of $\mathcal{A}$, and the spectral radius $1$ of $\mathcal{L}$ may be greater than one, and may not be equal to each other. In this case in principle it is possible that two different non-negative (normalised) eigenvectors of $\mathcal{A}$ corresponding to the spectral bound $0$ determine the same non-negative eigenvector of $\mathcal{L}$ corresponding to its spectral radius $1$, via formula \eqref{steadyconn2}. 
Hence we cannot conclude that the steady state problems in the two different formulations are equivalent. 
\end{enumerate} 
\begin{remark}
Note that, as we will see later in Section 5, we will establish the equi\-va\-lence of the steady state problems in the two different formulations, using a different approach. However, the approach and the results we established here are actually applicable to the nonlinear version of  model \eqref{distributed}. That is, when the vital rates depend on the total population size $P$. In that case the only (formal) difference is that at the steady state we deal with families of operators $\mathcal{A}_P$ and $\mathcal{L}_P$ parametrised by the total population size $P$, see e.g. \cite{AF,FGH} for more details. At the same time we note that the extension of the results of Section 5 to nonlinear (in particular quasilinear) models, promises to be a challenging problem.
\end{remark}

\section{On the equivalence of the time dependent problems}

Next we study the connection between the PDE and the delayed integral equation formulation of the linear distributed states at birth model, via establishing the connection between solutions of the two models. Let us start by recalling from \cite{D1} the formal definition of the initial value problem for the delay equation \eqref{inteq1} on the state space $L^1\left([-\Gamma,\infty);L^1(0,1)\right)$, with initial condition $\phi\in L^1\left([-\Gamma,0];L^1(0,1)\right)$ .
\begin{align}
b(\cdot,t)= & \int_0^1\frac{\beta(\cdot,y)}{\gamma(y)}\int_0^yb\left(x,t-\int_x^y\frac{1}{\gamma(r)}\,\ud r\right)
\exp\left\{-\int_x^y\frac{\mu(r)}{\gamma(r)}\,\ud r\right\}\,\ud x\,\ud y,\quad t>0, \nonumber \\
b(\cdot,t)= & \,\phi(\cdot,t),\quad t\in [-\Gamma,0]. \label{de}
\end{align}

Let us also recall from \cite{D1} the definition of a solution of the delay equation \eqref{de}. For a given
$\phi\in L^1\left([-\Gamma, 0];L^1(0,1)\right)$, $b\in L^1_{\text{loc}}\left([-\Gamma,\infty);L^1(0,1)\right)$ is a solution
of the delay equation \eqref{inteq1} if $b(\cdot,t)=\phi(\cdot,t)$ for $t\in[-\Gamma,0]$, and $b$ satisfies \eqref{de} for $t>0$.

We formulate our main results concerning the connection between solutions of the initial value problem for the
linear partial differential equation \eqref{distributed} and the delay equation \eqref{de} via two theorems. 
To establish the connection between \eqref{distributed} and \eqref{de} we consider an auxiliary problem, namely, the following inhomogeneous linear initial value problem.
\begin{equation}
\frac{ \ud p}{\ud t}=\mathcal{B}\, p+f,\quad t>0, \quad
p(0)=p^0,\label{inhomog}
\end{equation}
where $\mathcal{B}$ is defined in \eqref{Bdef}, and $f\in L^1\left((0,\infty);L^1(0,1)\right)$. Recall, for example from
\cite{P}, the definition of a mild solution $p$ of \eqref{inhomog}.
If $\mathcal{B}$ is the generator of a $C_0$ semigroup $\mathcal{T}$ on $L^1(0,1)$, then for $p^0\in L^1(0,1)$, the function $p$ given as
\begin{equation}
p(t)=\mathcal{T}(t)\, p^0+\int_0^t \mathcal{T}(t-\tau)f(\tau)\,\ud\tau,\quad 0\le t, \label{mild}
\end{equation}
is called the mild solution of \eqref{inhomog}.

Let us also recall a usual definition when dealing with delay equations. Namely, for a function $f \in
L^1([-\Gamma,\infty);L^1(0,1))$ and for $t \geq 0,$ we define $f_t \in L^1([-\Gamma,0];L^1(0,1))$ as $f_t(\tau)=f(t+\tau),\, \tau \in[-\Gamma,0].$ In particular, notice that $f_0=f|_{[-\Gamma,0]}.$

Now we also define a linear operator $K$ from $L^1([-\Gamma,0];L^1(0,1))$ to $L^1(0,1).$
\begin{equation}\label{K-opdef}
\left(Ku\right)(\cdot)=\frac{1}{\gamma(\cdot)}\int_0^\cdot
u\left(x,-\int_x^\cdot \frac{1}{\gamma(r)}\,\ud
r\right)\exp\left\{-\int_x^\cdot\frac{\mu(r)}{\gamma(r)}\,\ud
r\right\}\,\ud x.
\end{equation}
Our main tool to establish the connection between solutions of the delay and the partial differential equation is the following lemma. 
\begin{lemma}\label{newlemma}
Let us assume that $b \in L^1([-\Gamma,\infty),L^1(0,1)).$ Then the
function $t \rightarrow K\,b_t,\, t\in [0,\infty)$, is the mild solution of
the inhomogeneous linear initial value problem
\begin{equation}
\frac{ \ud p}{\ud t}=\mathcal{B}\, p+b,\quad t>0, \quad p(0)=K\,b_0.
\label{newhomog}
\end{equation}
\end{lemma}
{\bf Proof.} As we can see from \eqref{mild}, the mild solution depends on the semigroup $\mathcal{T}$ generated by $\mathcal{B}$. Hence we first need to determine explicitly the semigroup $\mathcal{T}$. To this end we consider the
initial value problem \eqref{inhomog} with $f\equiv 0$. In fact, to obtain the semigroup $\mathcal{T}$, we solve the corresponding homogeneous partial differential equation
\begin{equation}\label{homog-pde}
\begin{aligned}
p_t(s,t)+\left(\gamma(s)p(s,t)\right)_s &= -\mu(s)p(s,t),\quad s\in (0,1),  \\
\gamma(0)p(0,t)&=0, 
\end{aligned}
\end{equation} 
by integrating along characteristic curves. The characteristic curves are determined by the following system of ordinary
differential equations.
\begin{equation}\label{characteristics}
\dot{s}(t)=\gamma(s(t)),\quad
\dot{p}(t)=-(\mu(s(t))+\gamma'(s(t)))p(t),\quad s(0)=s^0,\quad p(0)=p^0\left(s^0\right).
\end{equation}
The solution of system \eqref{characteristics} is given by
\begin{equation}\label{characteristics2}
\left(S\left(t;s^0\right), p^0\left(s^0\right)\exp\left\{-\int_0^t\left(\mu\left(S(\tau;s^0)\right)+\gamma'\left(S(\tau;s^0)\right)\right)\,\ud\tau\right\}\right),
\end{equation}
where the function $S$ is implicitly defined via
\begin{equation}\label{S-def}
\int_{s^0}^{S\left(t;s^0\right)}\frac{1}{\gamma(r)}\,\ud r=t.
\end{equation}
Note that $S\left(t;s^0\right)$ is the size of individuals at time $t$, who were of size $s^0$ at time $0$. Hence the solution of the homogeneous partial differential equation \eqref{homog-pde} is given by
\begin{equation}\label{characteristic3}
p(s,t)=p^0(S^0(t;s))\exp\left\{-\int_0^t\left(\mu\left(S\left(\tau;S^0(t;s)\right)\right)
+\gamma'\left(S\left(\tau;S^0(t;s)\right)\right)\right)\,\ud\tau\right\},
\end{equation}
for $t\le\displaystyle\int_0^s\frac{1}{\gamma(r)}\,\ud r$, and $p(s,t)=0$
otherwise; where the function $S^0$ is implicitly defined via
\begin{equation*}
\int_{S^0(t;s)}^{s}\frac{1}{\gamma(r)}\,\ud r=t.
\end{equation*}
Note that $S^0(t;s)$ is the size of an individual at time $0$, who is of size $s$ at time $t$. Introducing the new integration variable $r=S(\tau;S^0(t;s))$ in \eqref{characteristic3} yields the following formula for the semigroup $\mathcal{T}$.
\begin{equation}\label{semigroup}
(\mathcal{T}(t)p)(s)=
\begin{Bmatrix}
0,\quad & t>\displaystyle\int_0^s\frac{1}{\gamma(r)}\,\ud r  \\
p\left(S^0(t;s)\right)\frac{\gamma\left(S^0(t;s)\right)}{\gamma(s)}\exp\left\{-\displaystyle\int_{S^0(t;s)}^s\frac{\mu(r)}{\gamma(r)}\,\ud
r\right\},\quad & t\le\displaystyle\int_0^s\frac{1}{\gamma(r)}\,\ud r
\end{Bmatrix}.
\end{equation}
Next we show that indeed $\mathcal{B}$ is the generator of the
semigroup $\mathcal{T}$ given in \eqref{semigroup}. To this end we prove that
for any $p\in D(\mathcal{B}),$ the right-hand derivative of $\mathcal{T}(t)\,p$ at $t=0$ is indeed $\mathcal{B}.$ That is,
\begin{equation}
\lim_{t\to 0^+}\displaystyle\int_0^1\left|\frac{(\mathcal{T}(t)p)(s)-p(s)}{t}-(\mathcal{B}p)(s)\right|\,\ud s=0,\quad p\in D(\mathcal{B}).
\end{equation}
On the one hand, utilising that the condition
$t>\displaystyle\int_0^s\frac{1}{\gamma(r)}\,\ud r$ is equivalent to $s<S(t;0),$ we have
\begin{align}
& \displaystyle\int_0^{S(t;0)}\left|\frac{(\mathcal{T}(t)p)(s)-p(s)}{t}-(\mathcal{B}p)(s) \right|\,\ud s=\int_0^{S(t;0)}\left|-\frac{p(s)}{t}-(\mathcal{B}p)(s)\right|\,\ud s \\
& \le \frac{1}{t}\int_0^{S(t;0)} |p(s)|\,\ud
s+\int_0^{S(t;0)}|(\mathcal{B}p)(s)|\,\ud s=
\frac{S(t;0)}{t}\left|p(\xi)\right|+\int_0^{S(t;0)}|(\mathcal{B}p)(s)|\,\ud
s,\label{diff-est}
\end{align}
for some $\xi\in [0,S(t;0)]$ (recall that $W^{1,1}$ functions are continuous). From \eqref{S-def} we have that $\frac{S(t;0)}{t}\le \displaystyle\max_{r\in [0,1]}\gamma(r)$, and since $p\in D(\mathcal{B})$, both terms in the right-hand side
of \eqref{diff-est} tend to $0$ as $t\to 0^+$.

On the other hand from \eqref{semigroup} we have
\begin{align}
& \frac{p\left(S^0(t;s)\right)\frac{\gamma\left(S^0(t;s)\right)}{\gamma(s)}\exp\left\{-\displaystyle\int_{S^0(t;s)}^s\frac{\mu(r)}{\gamma(r)}\,\ud r\right\}-p(s)}{t}-(\mathcal{B}p)(s) \nonumber \\
=& \frac{p\left(S^0(t;s)\right)\exp\left\{-\displaystyle\int_{S^0(t;s)}^s\frac{\mu(r)+\gamma'(r)}{\gamma(r)}\,\ud r\right\}-p(s)}{t}-(\mathcal{B}p)(s) \nonumber \\
= & \exp\left\{-\int_{S^0(t;s)}^s\frac{\mu(r)+\gamma'(r)}{\gamma(r)}\,\ud r\right\}\frac{p(S^0(t;s))-p(s)}{S^0(t;s)-s}\frac{S^0(t;s)-s}{t}+\gamma(s)p'(s) \nonumber \\
& +
p(s)\left(\frac{\exp\left\{-\int_{S^0(t;s)}^s\frac{\mu(r)+\gamma'(r)}{\gamma(r)}\,\ud
r\right\}-1}{t}+\mu(s)+\gamma'(s)\right)\to 0,
\end{align}
as $t\to 0^+$, for almost every $s\in\left[S^0(t;0),1\right]$. This is because $p\in D(\mathcal{B})$, and therefore the pointwise derivative of $p$ coincides with its distributional derivative almost everywhere (see for instance Theorem 1 in \cite[Sect. 1.1.3]{Maz}), and 
\begin{equation*}
\left.\frac{\partial}{\partial t}S^0(t;s)\right|_{t=0}=-\gamma(s).
\end{equation*} 
Hence we have
\begin{equation}
\displaystyle\int_{S(t;0)}^1\left|\frac{(\mathcal{T}(t)p)(s)-p(s)}{t}-(\mathcal{B}p)(s)
\right|\,\ud s\to 0,
\end{equation}
as $t\to 0^+$, by the Lebesgue dominated convergence theorem; in particular, since the difference quotients of a
$W^{1,1}$ function are bounded uniformly by the $L^1$ norm of its derivative.

We now compute the mild solution of \eqref{inhomog} by means of the variation of constants formula \eqref{mild}, with
inhomogeneity and initial condition given by \eqref{newhomog}.
Using the semigroup $\mathcal{T}$ given in \eqref{semigroup} we readily have, writing $b(s,\tau)$ for $b(\tau)(s),$
\begin{align}
& (\mathcal{T}(t)p^0)(s) \nonumber \\
& =\begin{Bmatrix}
0,\quad & t>\displaystyle\int_0^s\frac{1}{\gamma(r)}\,\ud r \\
\frac{1}{\gamma(s)}\displaystyle\int_0^{S^0(t;s)}b\left(\sigma,-\displaystyle\int_{\sigma}^{S^0(t;s)}\frac{1}{\gamma(r)}\,\ud
r\right)\exp\left\{-\displaystyle\int_\sigma^s\frac{\mu(r)}{\gamma(r)}\,\ud
r\right\}\,\ud\sigma,\quad & t\le \displaystyle\int_0^s\frac{1}{\gamma(r)}\,\ud r\label{sem}
\end{Bmatrix}.
\end{align}
Similarly, we have
\begin{align}
& \left(\mathcal{T}(t-\tau)b(\tau)\right)(s) \nonumber \\
& =\begin{Bmatrix}
0,\quad & t-\tau>\displaystyle\int_0^s\frac{1}{\gamma(r)}\,\ud r \\
b\left(S^0(t-\tau;s),\tau\right)\frac{\gamma\left(S^0(t-\tau;s)\right)}{\gamma(s)}\exp\left\{-\displaystyle\int_{S^0(t-\tau;s)}^s\frac{\mu(r)}{\gamma(r)}\,\ud
r\right\},\quad & t-\tau\le \displaystyle\int_0^s\frac{1}{\gamma(r)}\,\ud
r\label{var}
\end{Bmatrix}.
\end{align}
From \eqref{var} we have
\begin{align}
& \int_0^t\left(\mathcal{T}(t-\tau)b(\tau)\right)(s)\,\ud \tau \nonumber \\
= & \displaystyle\int\limits_{\displaystyle\max\left\{0,t-\int_0^s\frac{1}{\gamma(r)}\,\ud r\right\}}^{t} b\left(S^0(t-\tau;s),\tau\right)\frac{\gamma\left(S^0(t-\tau;s)\right)}{\gamma(s)}\exp\left\{-\int_{S^0(t-\tau;s)}^s\frac{\mu(r)}{\gamma(r)}\,\ud r\right\}\,\ud\tau \nonumber \\
= &
\frac{1}{\gamma(s)}\int\limits_{\displaystyle\max\left\{S^0(t;s),0\right\}}^{s} b\left(\sigma,t-\int_\sigma^s\frac{1}{\gamma(r)}\,\ud
r\right)\exp\left\{-\int_\sigma^s\frac{\mu(r)}{\gamma(r)}\,\ud
r\right\}\,\ud \sigma, \label{var2}
\end{align}
where to obtain the last equality we introduced the new variable
$\sigma:=S^0(t-\tau;s)$.

Adding \eqref{sem} and \eqref{var2} together, and using
\begin{equation*}
-\displaystyle\int_\sigma^{S^0(t;s)}\frac{1}{\gamma(r)}\,\ud r=t-\displaystyle\int_\sigma^s\frac{1}{\gamma(r)}\,\ud r,
\end{equation*} we have 
\begin{equation}
p(t):=\mathcal{T}(t)p^0+\int_0^t\mathcal{T}(t-\tau)b(\tau)\,\ud \tau=Kb_t,\quad t\in
[0,\infty). \label{loop}
\end{equation}
That is, the function $K\,b_t$ is the mild solution of the inhomogeneous initial value problem \eqref{newhomog}.\eofproof

Now we are in the position to formulate our first main result regarding the connection of solutions of the partial differential equation \eqref{distributed} and the delay equation \eqref{de}.
\begin{theorem}\label{theorem1}
For a given $\phi$, if $b$ is the unique non-negative solution of \eqref{de}, then the (unique) mild solution $p$ of
\eqref{inhomog} with
\begin{equation}
f=b|_{t\in (0,\infty)},\quad
p^0(\cdot)=\frac{1}{\gamma(\cdot)}\int_0^\cdot
b\left(x,-\int_x^\cdot\frac{1}{\gamma(r)}\,\ud r\right)
\exp\left\{-\int_x^\cdot\frac{\mu(r)}{\gamma(r)}\,\ud r\right\}\,\ud
x, \label{formula}
\end{equation}
i.e., the function given by formula \eqref{mild}, is also the (unique) mild (i.e. the semigroup) solution of \eqref{distributed} with initial condition $p^0$.
\end{theorem}
{\bf Proof.}
Using the definitions of the linear operators $\mathcal{C}$ in \eqref{Cdef}, and $K$ in \eqref{K-opdef}, 
for the solution $b(t):=b(\cdot,t)$ of \eqref{de} we have by hypothesis that $b(t)=\mathcal{C}K\,b_t$. 
On the other hand, since $p^0=K\,b_0,$ from Lemma \ref{newlemma}, it follows that $p(t)=K\,b_t.$
So,
\begin{align}
& \mathcal{T}(t)\,p^0+\int_0^t\mathcal{T}(t-\tau)\,\mathcal{C}(p(\tau))\,\ud \tau=\mathcal{T}(t)\,p^0+\int_0^t\mathcal{T}(t-\tau)\,\mathcal{C}\,K\,b_{\tau}\,\ud \tau \nonumber \\
&=\mathcal{T}(t)\,p^0+\int_0^t\mathcal{T}(t-\tau)\,b(\tau)\,\ud\tau=p(t),
\end{align}
where the last equality is just the definition of the mild solution of \eqref{inhomog}. Thus, $p(t)$ is the unique solution of the variation of constants equation corresponding to problem \eqref{distributed}. Hence $p(t)$ is the (semigroup) mild solution
of problem \eqref{distributed} (see also its abstract formulation \eqref{Cauchy-abstr}), see for example \cite[Sect.3.1]{P}.\eofproof

Our second result is formulated in the following theorem.
\begin{theorem}\label{theorem2}
Let us assume that $p^0\in L^1_+(0,1)$ is such that the functional equation
\begin{equation}\label{in-cond}
p^0(\cdot)=\frac{1}{\gamma(\cdot)}\int_0^\cdot
\phi\left(x,-\int_x^\cdot\frac{1}{\gamma(r)}\,\ud
r\right)\exp\left\{-\int_x^\cdot\frac{\mu(r)}{\gamma(r)}\,\ud
r\right\}\,\ud x,
\end{equation}
has a non-negative solution $\phi\in
L^1\left((0,1)\times[-\Gamma,0]\right)$. 
Then, the unique non-negative mild solution $p$ of \eqref{distributed} with
$p(\cdot,0)=p^0(\cdot)$ determines a non-negative solution
$b$ of \eqref{de}, given by
\begin{equation}\label{fixed-point}
b(\cdot,t)=\int_0^1\beta(\cdot,y)p(y,t)\,\ud y,\,\quad t>0;\quad\quad b(\cdot,t)=\phi(\cdot,t),\quad t\in [-\Gamma,0].
\end{equation}
\end{theorem}
{\bf Proof.} The mild solution $p$ of \eqref{distributed} satisfies the variation of constants equation
\begin{equation}
p(t)=\mathcal{T}(t)\,p^0+\int_0^t\mathcal{T}(t-\tau)\,\left(\mathcal{C}p\right)(\tau)\,\ud\tau,
\end{equation}
hence it is also a mild solution of the inhomogeneous equation
\eqref{inhomog} (with the same initial condition), with
$f=\mathcal{C}\, p$. Now we define $b$ as in \eqref{fixed-point}, i.e., $b(t)=\mathcal{C}(p(t))$ for $t>0,$ and $b(t) = \phi(t)$ for $t \in [-\Gamma,0].$ It turns out then, that $p$ solves \eqref{newhomog} in the mild sense (note that $p^0 = K\phi = K\,b_0$ by \eqref{in-cond}). Hence, by Lemma \ref{newlemma},
\begin{equation}\label{kc-eq}
p(t)=K\, b_t, \quad t \geq 0,
\end{equation}
holds. From equation \eqref{kc-eq} we have now, for $t > 0,$
\begin{equation}
b(t)=\mathcal{C}(p(t))=\mathcal{C}(K\, b_t),
\end{equation}
and the proof is completed.\eofproof

The previous two theorems allow us to draw the following somewhat interesting conclusion. 
\begin{corollary}
Assume that the delay equation \eqref{de} has two solutions $b_1$ and $b_2$, such that $\phi_1=b_1|_{t\in [\Gamma,0]}$ and $\phi_2=b_2|_{t\in [\Gamma,0]}$ determine the same function $p^0$ via the right-hand side of equation \eqref{in-cond}. Then $b_1(\cdot,t)\equiv b_2(\cdot,t),\,\forall t>0$ holds.
\end{corollary}

We note that it can be shown that for any non-negative integrable function $\phi$, the right-hand side of \eqref{in-cond} determines a non-negative integrable function $p^0$.

Intuitively it is also clear, that for any given initial size-distribution $p^0$ we can calculate what was the size of an individual, which at time $0$ is of size $y$, at any given time $t\in [-\Gamma,0]$, (but note that we need to make
sure that the size remains positive). That is, it is shown that for any initial density $p^0\in L^1_+(0,1)$, equation \eqref{in-cond} admits a non-negative solution $\phi\in L^1\left((0,1)\times[-\Gamma,0]\right)$. We show this here for the special case of $\gamma\equiv 1$ and $\mu\equiv 0$. For any given size-distribution $p^0\in L^1_+(0,1)$ let us choose an arbitrary age-size-distribution $P\in L^1_+\left((0,1)\times (0,1)\right)$, such that for any size $y\in (0,1)$ we have
\begin{equation}\label{agesize}
p^0(y)=\int_0^1P(y,a)\,\ud a.
\end{equation}
For example, we can choose 
\begin{equation*}
P(y,a)=\frac{p^0(y)}{y}\chi_{[0,y]}(a),\quad y,a\in (0,1).
\end{equation*} 
Further, we define $\phi\in L^1_+\left((0,1)\times(-1,0)\right)$ as
\begin{equation}
\phi(s,t):=\begin{Bmatrix}
P(s-t,-t) & \text{if}\quad s\le t+1 \\
0 & \text{if}\quad s>t+1 
\end{Bmatrix}.
\end{equation}
We then have 
\begin{equation}\label{agesize2}
\int_0^y\phi(x,x-y)\,\ud x =\int_0^yP(y,y-x)\,\ud x= \int_0^yP(y,a)\,\ud a=\int_0^1 P(y,a)\,\ud a=p^0(y),
\end{equation}
that is \eqref{in-cond} holds. Also we have
\begin{align}
& \int_{-1}^0\int_0^1\phi(s,t)\,\ud s\,\ud t \nonumber \\ 
&=\int_{-1}^0\int_0^{t+1}P(s-t,-t)\,\ud s\,\ud t=\int_0^1\int_0^yP(y,a)\,\ud a\,\ud y=\int_0^1p^0(y)\,\ud y<\infty,\label{agesize3}
\end{align}
hence $\phi\in  L^1_+\left((0,1)\times (-1,0)\right)$ indeed, as required. The existence of an appropriate $\phi$ for a given $p^0$ in the general case (i.e. when $\gamma \not\equiv 1$ and $\mu\not\equiv 0$) is shown similarly, but the formulas are much more cumbersome.

\section{General results}

Motivated by the results of the previous sections, in this section we establish the equivalence between a general class of PDE's and their corresponding delay (renewal) equation formulations. As we will see, the PDE model \eqref{distributed} discussed in the previous sections fits into the general framework we present here. However, in contrast to the previous section, we take a slightly different point of view. To be able to work with greater generality, our starting point is the semigroup $\mathcal{T}_0$ generated by an operator $\mathcal{B}$, which is assumed to describe individual development (e.g. individual growth) and survival. Note that the operator defined in \eqref{Bdef} is a particular example of such an operator $\mathcal{B}$ (that is why purposefully we use the same notation throughout the section). Here, to avoid cumbersome formulas and calculations, we assume that time is scaled such that the maximal age an individual may attain equals $1$. This then amounts to that the semigroup $\mathcal{T}_0$ is nilpotent, in particular $\mathcal{T}_0(t)\,u=0$ holds, for $t\ge 1$, for any $u\in\mathcal{X}$. We assume that in general reproduction (or recruitment of individuals) is described by a bounded linear operator $\mathcal{C}\,:\,\mathcal{X}\to\mathcal{X}$. Again, the operator defined in \eqref{Cdef} is a particular example of such a general recruitment operator $\mathcal{C}$.

Once the semigroup $\mathcal{T}_0$ is known (given), the semigroup $\mathcal{T}$ generated by $\mathcal{A}=\mathcal{B}+\mathcal{C}$ can be constructed by solving pointwise the variation of constants equation
\begin{equation}
\mathcal{T}(t)=\mathcal{T}_0(t)+\int_0^t\mathcal{T}_0(t-s)\,\mathcal{C}\,\mathcal{T}(s)\,\ud s. \label{DE1}
\end{equation} 
Applying $\mathcal{C}$ to \eqref{DE1} we obtain
\begin{equation}
\mathcal{L}(t)=\mathcal{L}_0(t)+\int_0^t\mathcal{L}_0(t-s)\,\mathcal{L}(s)\,\ud s,\label{DE2}
\end{equation}
where we defined
\begin{equation*}
\mathcal{L}_0(t):=\mathcal{C}\,\mathcal{T}_0(t),\quad \text{and}\quad \mathcal{L}(t):=\mathcal{C}\,\mathcal{T}(t).
\end{equation*}
The solution of equation \eqref{DE2} is obtained by generation expansion as
\begin{equation*}
\mathcal{L}=\displaystyle\sum_{n=1}^{\infty}\mathcal{L}^{n\otimes},
\end{equation*}
see Section 2 of \cite{DGM1} for more details. Once $\mathcal{L}$ is known, equation \eqref{DE1} becomes an explicit formula for the semigroup $\mathcal{T}$,
\begin{equation}
\mathcal{T}(t)=\mathcal{T}_0(t)+\int_0^t\mathcal{T}_0(t-s)\mathcal{L}(s)\,\ud s.\label{DE3}
\end{equation}
\begin{lemma}\label{Lemma1}
We have
\begin{equation}
(\lambda\,\mathcal{I}-\mathcal{A})^{-1}=(\lambda\,\mathcal{I}-\mathcal{B})^{-1}\left(\mathcal{I}-\hat{\mathcal{L}_0}(\lambda)\right)^{-1},\label{DE4}
\end{equation}
where $\hat{\mathcal{L}_0}(\lambda)$ stands for the Laplace transform of $\mathcal{L}(t)$.
\end{lemma}
{\bf Proof.} Taking the Laplace transform of \eqref{DE2} we have
\begin{equation}
\hat{\mathcal{L}}(\lambda)=\left(\mathcal{I}-\hat{\mathcal{L}_0}(\lambda)\right)^{-1}\hat{\mathcal{L}_0}(\lambda).\label{DE5}
\end{equation}
Laplace transforming \eqref{DE3}, and using that the Laplace transform of a strongly conti\-nuous semigroup coincides with the resolvent operator of its generator (see for example \cite[Sect.II]{NAG}), yields
\begin{equation}
(\lambda\,\mathcal{I}-\mathcal{A})^{-1}=(\lambda\,\mathcal{I}-\mathcal{B})^{-1}\left(\mathcal{I}+\hat{\mathcal{L}}(\lambda)\right).\label{DE6}
\end{equation}
Combining formulas \eqref{DE5} and \eqref{DE6} yields \eqref{DE4}. \eofproof
\begin{corollary}
\begin{equation}
\sigma(\mathcal{A})=\left\{\lambda\in\mathbb{C}\,|\,1\in\sigma\left(\hat{\mathcal{L}_0}(\lambda)\right)\right\}.
\end{equation}
\end{corollary}
{\bf Proof.} Since $\mathcal{T}_0$ is nilpotent, its resolvent $(\lambda\,\mathcal{I}-\mathcal{B})^{-1}$ is analytic on $\mathbb{C}$. \eofproof

Let us now turn to the delay formulation of the abstract Cauchy problem 
\begin{equation}\label{Cauchy}
\frac{\ud p}{\ud t}=\mathcal{A}\,p,\quad p(0)=p^0.
\end{equation} 
By this we mean the renewal equation (RE)
\begin{equation}
b(t)=\int_0^1\mathcal{L}_0(a)\,b(t-a)\,\ud a, \label{DE7}
\end{equation}
together with the initial-history condition
\begin{equation}
b(\theta)=\phi(\theta),\quad -1\le\theta\le 0,\label{DE8}
\end{equation}
with $\phi\in L^1\left([-1,0];\mathcal{X}\right)=:\mathcal{Y}$. 
Note that the delay equation \eqref{de} becomes a special case of \eqref{DE7}-\eqref{DE8}. That is, \eqref{DE7}-\eqref{DE8} coincides with \eqref{de}, when the operators $\mathcal{B}$ and $\mathcal{C}$ are defined as in \eqref{Bdef}, \eqref{Cdef}.

The solution of \eqref{DE7}-\eqref{DE8} is given by
\begin{equation}
b(t)=f(t)+\int_0^t\mathcal{L}(a)f(t-a)\,\ud a,\quad t>0,\label{DE9}
\end{equation}
where
\begin{equation}
f(t)=f(t,\phi):=\int_t^1\mathcal{L}_0(a)\phi(t-a)\,\ud a,\label{DE10}
\end{equation}
with the understanding that the integral above in \eqref{DE10} equals zero for $t\ge 1$. 

The definition
\begin{equation}
\left(\mathcal{S}(t)\,\phi\right)(\theta)=b(t+\theta;\phi)=b_t(\theta;\phi)\label{DE11}
\end{equation}
yields a strongly continuous semigroup on $\mathcal{Y}$.
A straightforward Laplace transform calculation shows that
\begin{equation}
(\lambda\mathcal{I}-\tilde{\mathcal{A}})^{-1}\phi(\theta)=e^{\lambda\theta}c+\int_{\theta}^0e^{\lambda(\theta-\sigma)}\phi(\sigma)\,\ud \sigma,\label{DE13}
\end{equation}
where $\tilde{A}$ is the generator of $\mathcal{S}$, and
\begin{equation}
c=\left(\mathcal{I}-\hat{\mathcal{L}_0}(\lambda)\right)^{-1}\int_0^1\mathcal{L}_0(a)\int_{-a}^0e^{-\lambda(a+\sigma)}\phi(\sigma)\,\ud \sigma\,\ud a.\label{DE14}
\end{equation}
\eqref{DE13}-\eqref{DE14} then shows that the spectra of $\mathcal{A}$ and $\tilde{\mathcal{A}}$ coincide, which already indicates a strong connection between the semigroups $\mathcal{S}$ and $\mathcal{T}$. We now further explore this connection. 

We define the linear operator $K\,:\,\mathcal{Y}\to\mathcal{X}$ as
\begin{equation}
K\,\phi=\int_{-1}^0\mathcal{T}_0(-s)\,\phi(s)\,\ud s=\int_0^1\mathcal{T}_0(a)\,\phi(-a)\,\ud a.\label{DE15}
\end{equation}
Note that the operator $K$ defined in \eqref{K-opdef} is a special case of $K$ defined above in \eqref{DE15}. 
\begin{lemma}\label{lemma-DE}
\begin{equation}
b(t;\phi)=\mathcal{C}\,\mathcal{T}(t)\,(K\,\phi).\label{DE16}
\end{equation}
\end{lemma}
{\bf Proof.} Recall \eqref{DE10}, and recall that $\mathcal{T}_0=0$ for $t\ge 1$. It follows that
\begin{align}
f(t;\phi)& =\int_t^{1+t}\mathcal{L}_0(a)\,\phi(t-a)\,\ud a=\int_{-1}^0\mathcal{L}_0(t-\sigma)\,\phi(\sigma)\,\ud \sigma \nonumber \\
& =\mathcal{L}_0(t)\int_{-1}^0\mathcal{T}_0(-\sigma)\,\phi(\sigma)\,\ud\sigma=\mathcal{L}_0(t)\,(K\,\phi),
\end{align}
where we used that $\mathcal{L}_0(t-\sigma)=\mathcal{C}\,\mathcal{T}_0(t-\sigma)=\mathcal{C}\,\mathcal{T}_0(t)\,\mathcal{T}_0(-\sigma)=\mathcal{L}_0(t)\,\mathcal{T}_0(-\sigma)$. 

Inserting the expression for $f$ into \eqref{DE9} we obtain
\begin{equation}
b(t;\phi)=\mathcal{L}_0(t)\,K\,\phi+\int_0^t\mathcal{L}(a)\mathcal{L}_0(t-a)\,(K\,\phi)\,\ud a.\label{DE17}
\end{equation}
Similarly to \eqref{DE1}-\eqref{DE2}, recalling $\mathcal{T}(t)=\mathcal{T}_0(t)+\int_0^t\mathcal{T}(s)\,\mathcal{C}\,\mathcal{T}_0(t-s)\,\ud s$, and its counterpart $\mathcal{L}(t)=\mathcal{L}_0(t)+\int_0^t\mathcal{L}(s)\,\mathcal{L}_0(t-s)\,\ud s$ , we have
\begin{equation}
b(t;\phi)=\mathcal{L}(t)\,K\,\phi=\mathcal{C}\,\mathcal{T}(t)\,(K\,\phi).
\end{equation}
\eofproof

Note the similarity between Lemma \ref{newlemma} and Lemma \ref{lemma-DE} above. The connection between the abstract Cauchy problem \eqref{Cauchy} and the renewal equation \eqref{DE7}-\eqref{DE8} is established in the following theorem. 
\begin{theorem}\label{equivalence}
\begin{equation}
\mathcal{T}(t)\,(K\,\phi)=K\,\mathcal{S}(t)\,\phi.
\end{equation}
\end{theorem}
{\bf Proof.} We have
\begin{equation*}
\mathcal{T}(t)\,(K\,\phi)=\mathcal{T}_0(t)\, (K\,\phi)+\int_0^t\mathcal{T}_0(t-s)\,\mathcal{L}(s)\,(K\,\phi)\,\ud s.
\end{equation*}
Since
\begin{equation}
\mathcal{T}_0(t)\,(K\,\phi)=\int_{-1}^0\mathcal{T}_0(t-s)\,\phi(s)\,\ud s,\quad \text{and}\quad \mathcal{L}(s)\,(K\,\phi)=b(s),
\end{equation}
we have
\begin{align}
\mathcal{T}(t)\,(K\,\phi) & = \int_{-1}^t\mathcal{T}_0(t-s)\,b(s)\,\ud s=\int_{t-1}^t\mathcal{T}_0(t-s)\,b(s)\,\ud s \nonumber \\
& = \int_{-1}^0\mathcal{T}_0(-\sigma)\,b(t+\sigma)\,\ud \sigma=K\,b_t=K\,\mathcal{S}(t)\,\phi.
\end{align} \eofproof

Again, note the similarity between Theorems \ref{theorem1}-\ref{theorem2} and Theorem \ref{equivalence}.

In the lemma below $\mathcal{R}(\mathcal{O})$ stands for the range of an operator $\mathcal{O}$.
\begin{lemma}\label{Lemma5}
For $t\ge 1$ the inclusion $\mathcal{R}(\mathcal{T}(t))\subset\mathcal{R}(K)$ holds.
\end{lemma}
{\bf Proof.} For $t\ge 1$ the identity \eqref{DE3} reduces to
\begin{equation*}
\mathcal{T}(t)=\int_0^t\mathcal{T}_0(t-s)\,\mathcal{L}(s)\,\ud s,
\end{equation*}
and
\begin{align*}
\int_0^t\mathcal{T}_0(t-s)\,\mathcal{L}(s)\,\ud s &=\int_{t-1}^t\mathcal{T}_0(t-s)\,\mathcal{L}(s)\,\ud s=\int_{-1}^0\mathcal{T}_0(-\sigma)\,\mathcal{L}(t+\sigma)\,\ud \sigma  \\
\Rightarrow\, \forall\,\psi\in\mathcal{X},\quad \mathcal{T}(t)\,\psi &=K\left(\underbrace{\theta\mapsto\mathcal{L}(t+\theta)\,\psi}_{\phi}\right),\quad \text{for}\quad t\ge 1.
\end{align*}
\eofproof

\begin{remark}
In view of Lemma \ref{Lemma5}, Theorem \ref{equivalence} tells us that the large time behaviour of $\mathcal{T}$ can be described in terms of the large time behaviour of $\mathcal{S}$.
\end{remark}

The equivalence of the steady state problems in the two different formulations is established in the following theorem. Below  $\overline{p}$ and $\overline{\phi}$ stand for time-independent solutions of \eqref{Cauchy} and \eqref{DE7}-\eqref{DE8}, respectively.
\begin{theorem}\label{ss-equivalence} 
If $\mathcal{S}(t)\,\overline{\phi}=\overline{\phi}$, then $\mathcal{T}(t)\,\overline{p}=\overline{p}$, where $\overline{p}=K\,\overline{\phi}$. On the other hand, if $\mathcal{T}(t)\,\overline{p}=\overline{p}$, then $\mathcal{S}(t)\,\overline{\phi}=\overline{\phi}$, with $\overline{\phi}=\mathcal{C}\,\overline{p}$, i.e. $\overline{\phi}(\theta)=\mathcal{C}\,\overline{p}$, $-1\le\theta\le 0$.
\end{theorem}
{\bf Proof.} First note that if $\overline{\phi}$ is such that $\mathcal{S}(t)\,\overline{\phi}=\overline{\phi}$ holds, then
\begin{equation*}
\mathcal{T}(t)\,\overline{p}=\mathcal{T}(t)\,K\,\overline{\phi}=K\,\mathcal{S}(t)\,\overline{\phi}=K\,\overline{\phi}=\overline{p},
\end{equation*} 
shows that $\overline{p}=K\,\overline{\phi}$ is a steady state of \eqref{Cauchy}.

On the other hand, for a time-independent solution $\overline{p}$ of \eqref{Cauchy}, for $t\ge 1$ we have
\begin{align*}
\overline{p} &=\mathcal{T}(t)\,\overline{p}=\mathcal{T}_0(t)\,\overline{p}+\int_0^t\mathcal{T}_0(t-s)\,\mathcal{C}\,\overline{p}\,\ud s  \\
&=\int_{t-1}^t\mathcal{T}_0(t-s)\,\ud s\,\mathcal{C}\,\overline{p}=\int_0^1\mathcal{T}_0(a)\,\ud a\,\mathcal{C}\,\overline{p}.
\end{align*}
Applying $\mathcal{C}$ to the identity above we obtain
\begin{equation}
\mathcal{C}\,\overline{p}=\mathcal{C}\,\overline{p}\,\int_0^1\mathcal{L}_0(a)\,\ud a,
\end{equation} 
which means that $\mathcal{C}\,\overline{p}$ is a constant solution of the renewal equation \eqref{DE7}-\eqref{DE8}, and the proof is completed. \eofproof

\section{Concluding remarks}
In this work first we considered a structured population model with distributed states at birth. The model can naturally be formulated as a first order partial differential equation on the appropriate state space of Lebesgue integrable functions. At the same time we note that the  model belongs to a general class of physiologically structured population mo\-dels discussed for example in \cite{DGM1}. Our main goal was to investigate how one can recover the corresponding delay (renewal) equation formulation of the model, similarly to the Volterra integral equation formulation of a model with single state at birth. In\-te\-restingly enough, the connection between the two formulations can be motivated (at least in our opinion) by a fairly recent spectral theoretic result in \cite{HT}, due to H. Thieme. In particular, as we have seen in Section 2, Theorem \ref{spectralth} provides a rigorous  connection between two spectral problems arising from the different formulations at the steady state, at least when some appropriately defined operators are irreducible. In general, irreducibility plays an important role when transforming steady state problems (also of nonlinear models) into spectral problems. It guarantees the existence and  uniqueness of a (normalised) strictly positive eigenvector corresponding to the spectral bound (or radius) of an operator. Partly for this reason, we explored the connection between the irreducibility of the semigroup and the integral operator, arising from the two different formulations of the steady state problem. We showed that irreducibility of the integral operator implies irreducibility of the semigroup, but not vice versa. We also note that Theorem \ref{spectralth} already shows that the PDE and delay formulations are ``asymptotically equivalent''. That is, in the generic case (when the spectral bound of $\mathcal{A}$ is not zero) solutions of the PDE and the delay equation grow or decay simultaneously.
 
In Section 4 we established the connection between the two time-dependent problems via Theorems \ref{theorem1} and \ref{theorem2}. It turns out, that for any solution of the partial differential equation there exists a corresponding solution of the delayed integral equation and vice versa. As far as we know the results presented here are the first ones comparing  solutions of the delay and PDE formulations of a physiologically structured population model with distributed states at birth.
We also note that the results established in Section 4 can readily be extended (at least) for semilinear equations, when solutions of the PDE can be obtained using the method of characteristics. One anticipates that the results carry forward from the linear case, but the calculations and formulas will be much more cumbersome. We anticipate though, that a similar comparison of other classes of models, in particular quasi-linear ones, may reveal differences. A natural candidate would be a quasilinear hierarchic (infinite dimensional environment/interaction variable) size-structured population model with single state at birth. It was shown in \cite{AD,AI} that this model is not well-posed on the biologically natural state space  $L^1$, if the growth rate is not a monotone function of the environment. Spe\-ci\-fi\-cally, in \cite{AI} it was shown that finite time blow up of the density function is possible, and consequently, one has to work with measures to study existence of solutions. At the same time note that models of this type (for example the one studied in \cite{AI}) can be formulated and studied in the delay equation framework, too.  In particular, an interesting cannibalistic model with infinite dimensional environment was shown to be well-posed in \cite{D1} in the delay equation framework, under the assumptions that the vital rates are continuously differentiable. 

Motivated by the results obtained in Section 4 for the distributed states at birth model \eqref{distributed}, we established some very general results in Section 5. Our viewpoint (and starting point) in Section 5 was also slightly different to that in Section 4. In particular, to allow greater generality, and at the same time to avoid the delicate technical calculations of Section 4, our starting point in Section 5 is the semigroup governing a PDE model. This then allowed us to establish results in Section 5, which directly apply to the distributed states at birth model \eqref{distributed}. In particular, Theorem \ref{ss-equivalence} establishes the equivalence of the steady state problems in the two different formulations; while Theorem \ref{equivalence} establishes the rigorous connection between time-dependent solutions  of the two models. 

However, to extend the results of Section 5 for nonlinear models, promises to be rather challenging. There are also natural, biologically motivated examples of recruitment processes, which lead to an unbounded operator $\mathcal{C}$, posing additional difficulties.  We aim to work on these problems in the future. 

One may also ask the very na\-tu\-ral question what happens in a perhaps simpler single state at birth model. (Note that when we say "simpler" we mean from the modelling/biological point of view.) Indeed it will be very interesting to explore whether Theorem \ref{spectralth} could be generalised, and one could use a similar operator splitting in case of the single state at birth model, by considering the recruitment at state $0$ as a boundary perturbation. In fact there are numerous results, see e.g. \cite{D-Sch,NAG,Gre}, on boundary perturbations of strongly continuous semigroups. The main idea of the Desch-Schappacher boundary perturbation theory is to lift the problem to the extrapolated space where one defines an additive, bounded and positive perturbation. Then, the boundary perturbation, i.e. the recruitment operator, can be recovered as the part of the generator of the semigroup (defined on the extrapolated space) in the original state space.  We may anticipate that one will then apply a similar result to Theorem \ref{spectralth}, for the generator of the extrapolated semigroup. Since the extrapolated semigroup and the original semigroup are similar, their spectrum, spectral bound and growth bound coincide. We aim to elaborate  the details of this in the future.

\section*{Acknowledgments}
\`{A}. Calsina was partially supported by the research projects 2009SGR-345 and DGI MTM2011-27739-C04-02.
J.\,Z. Farkas was supported by the research project DGI MTM2011-27739-C04-02, while visiting the Universitat Aut\`{o}noma de Barcelona, and by a personal research grant from The Carnegie Trust for the Universities of Scotland.

\end{document}